\documentclass[11pt]{article}
\usepackage{fullpage,amsmath,graphicx,amsthm,enumitem}

\newtheorem{theorem}{Theorem}[section]
\newtheorem{lemma}[theorem]{Lemma}
\newtheorem{corollary}[theorem]{Corollary}
\newtheorem{prop}[theorem]{Proposition}
\numberwithin{figure}{section}

\newcommand{\de}{\backslash}

\begin{document}

\title{\bf Graphs without large $K_{2,n}$-minors}

\author{Guoli Ding
\\ Department of Mathematics, Louisiana State University, Baton Rouge, USA}

\date{\today}

\maketitle

\begin{abstract}
The purpose of this paper is to characterize graphs that do not have a large $K_{2,n}$-minor. As corollaries, it is proved that, for any given positive integer $n$, every sufficiently large 3-connected graph with minimum degree at least six, every 4-connected graph with a vertex of sufficiently high degree, and every sufficiently large 5-connected graph must have a $K_{2,n}$-minor.
\end{abstract}

\section{Introduction}

All graphs considered in this paper are simple, as loops and parallel edges do not affect if or not a graph has a $K_{2,n}$-minor. In particular, a minor should be interpreted as a simple minor. 

To explain the kind of theorems that we are going to prove, we first consider $K_{1,n}$, instead of $K_{2,n}$. This is because the structure is simpler in this case and thus it is easier to make a point. For each positive integer $n$, let ${\cal S}_n$ be the class of graphs that do not contain $K_{1,n}$ as a minor. Then for any graph $G$, it is obvious that $G\in{\cal S}_1$ if and only if each component of $G$ is $K_1$; $G\in{\cal S}_2$ if and only if each component of $G$ is either $K_1$ or $K_2$; and $G\in{\cal S}_3$ if and only if each component of $G$ is either a path or a cycle. It is also easy to characterize ${\cal S}_4$ as follows.

\begin{prop}\label{prop:s4}
A graph is in ${\cal S}_4$ if and only if each of its components is a minor of one of the following two graphs, where the thin edges could be subdivided any number of times.
\end{prop}

\begin{figure}[ht]
\centerline{\includegraphics[scale=0.54]{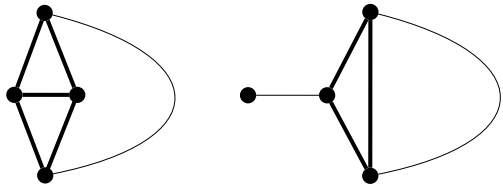}}
\caption{Maximal connected graphs that do not contain a $K_{1,4}$-minor}
\end{figure}

With some extra efforts, one can also characterize ${\cal S}_5$ in a similar way, where the number of maximal graphs increases from two to eighteen. These results suggest that other $\mathcal S_n$ could admit a similar characterization and it turns out that this is true. 

\begin{prop} \label{prop:sn}
There exists a function $s(n)$ such that every component of a graph in $\mathcal S_n$ is a subdivision of a connected graph on at most $s(n)$ vertices.
\end{prop}

This is a reformulation of a result of Robertson and Seymour \cite{RS}. It implies that every ${\cal S}_n$ can be characterized in a way similar to Proposition \ref{prop:s4}. However, it seems quite unlikely that we can obtain, like in Proposition \ref{prop:s4}, an explicit list of maximal graphs for a general $n$. Therefore, in a sense, Proposition \ref{prop:sn} is the best we can do to characterize a general ${\cal S}_n$. 

Strictly speaking, Proposition \ref{prop:sn} is not a characterization of ${\cal S}_n$. However, it does present a structure that every graphs in ${\cal S}_n$ must have. It is important to point out that this structure has the following nice property. Let ${\cal T}_n$ be the class of graphs for which each component is a subdivision of a graph on at most $n$ vertices. Then Proposition \ref{prop:sn} can be restated as: there exists a function $s(n)$ such that ${\cal S}_n \subseteq {\cal T}_{s(n)}$ holds for all $n$. On the other hand, it is not difficult to show that ${\cal T}_n\subseteq {\cal S}_{n^2}$ holds for all $n$. In other words, every graph in ${\cal S}_n$ must be in ${\cal T}_{s(n)}$ and, although ${\cal T}_{s(n)}$ may contain graphs that do have a $K_{1,n}$-minor, but no graph in ${\cal T}_{s(n)}$ has a $K_{1,m}$-minor when $m$ is sufficiently large. From this point of view, Proposition \ref{prop:sn} can be considered as an approximate characterization of every ${\cal S}_n$.

There are certainly ways to formulate Proposition \ref{prop:sn} such that it is a characterization, in the ordinary sense, of a graph property.

\begin{prop} 
A class of connected graphs is contained in some ${\cal S}_n$ if and only if it is contained in some ${\cal T}_m$.
\end{prop}

The next is another equivalent formulation of Proposition \ref{prop:sn}. Let us call a class of connected graphs {\it minor-closed} if every connected minor of a member is also a member. For each positive integer $r$, let ${\cal K}_r=\{K_{r,s}:s\ge r\}$.

\begin{prop} \label{prop:sniff} 
Let ${\cal G}$ be a minor-closed class of connected graphs. Then ${\cal K}_1 \not \subseteq {\cal G}$ if and only if ${\cal G} \subseteq {\cal T}_m$ for some $m$.\end{prop}

Very often, Proposition \ref{prop:sniff} is stated very loosely as: {\it a connected graph does not have a large $K_{1,n}$-minor if and only if it is obtained from a graph of bounded size by subdividing edges.}

The goal of this paper is to prove, for $K_{2,n}$, a theorem analogous to propositions \ref{prop:sn}-\ref{prop:sniff}. We will show that, loosely speaking, a 2-connected graph does not have a large $K_{2,n}$-minor if and only if it is obtained by taking ``2-sums", for a bounded number of iterations, of two types of graphs: one type is a generalization of outerplanar graphs and the other type is a kind of  ``subdivision" of graphs of bounded size. To state our theorem more precisely, we need to introduce some definitions. First, we define the two types of graphs from which we built all graphs that do not have a large $K_{2,n}$-minor. 

The first type is a generalization of outerplanar graphs. Let $G$ be a graph with a specified Hamiltonian cycle $C$, which will be called the {\it reference cycle}. Other edges of $G$ are called {\it chords} and two non-incident chords $ab$ and $cd$ {\it cross} if the four vertices appear on the cycle $C$ in the order $a,c,b,d,a$. We call $G$ {\it type-I} if every chord crosses at most one other chord and, in addition, if two chords $ab$ and $cd$ do cross, then either both $ac,bd$ are edges of $C$ or both $ad,bc$ are edges of $C$. Clearly, all 2-connected outerplanar graphs are of type-I. It is well known that outerplanar graphs do not have a $K_{2,3}$-minor. For type-I graphs, they may have a $K_{2,3}$-minor, but we will show, in Section 6, that they do not have $K_{2,5}$-minors. The class of all type-I graphs is denoted by $\cal P$.

Graphs of the second type are constructed from graphs of bounded size. To define these graphs we begin with two special classes of graphs. Let $H$ be a type-I graph, $C$ be the reference cycle of $H$, and  $ab,cd$ be two distinct edges of $C$. Suppose all chords of $C$ are between the two paths of $C\de\{ab,cd\}$. If $ab$ and $cd$ have a common end, say $a=d$, then $H$ is called a {\it fan} with {\it corners} $a,b,c$. Vertex $a$ will be called the {\it center} of the fan, and the number of chords will be called the {\it length} of the fan. If $ab$ and $cd$ have no common ends, then for any subset $F$ of $\{ab,cd\}$, $H\de F$ is called a {\it strip} with {\it corners} $a,b,c,d$, provided that the minimum degree of $H\de F$ is at least two. The {\it length} of a strip is the maximum number of chords that are pairwise noncrossing and pairwise nonadjacent. 

Let $G$ by a graph. Then {\it adding} a fan or a strip to $G$ means to identify the corners of a fan or a strip (which is disjoint from $G$) with distinct vertices of $G$. An {\it augmentation} of a graph $G$ is obtained by adding disjoint fans and strips to $G$ such that if two corners are identified with the same vertex of $G$ then one of them is the center of a fan and the other is either a center of the fan or a corner of a strip. For each positive integer $m$, let ${\cal B}_m$ be the class of graphs on at most $m$ vertices and let ${\cal A}_m$ be the class of all augmentations of graphs in ${\cal B}_m$. These are our second type of graphs. We will show, in Section 6, that if $n$ is sufficiently large, then no graph in ${\cal A}_m$ can have a $K_{2,n}$-minor.

Next, we describe an operation by which we can construct from $\cal P$ and $\mathcal A_m$ all graphs that do not have a large $K_{2,n}$-minor.

Let $G_1$ and $G_2$ be two disjoint graphs. For $i=1,2$, let $z_i$ be a vertex of $G_i$ such that it is incident with precisely two edges $x_iz_i$ and $y_iz_i$; let $G_i'=G_i\de z_i$ if $x_iy_i\not\in E(G_i)$ and $G_i'=G_i\de z_i\de e_i$ if $e_i=x_iy_i\in E(G_i)$. Then a 2-{\it sum} of $G_1$ and $G_2$, over $z_1$ and $z_2$, is a graph obtained from $G_1'$ and $G_2'$ by first identifying $x_1$ with $x_2$, $y_1$ with $y_2$, and then, in case $x_iy_i\in E(G_i)$ holds for at least one $i$, adding an edge between the two new vertices, which are the {\it joins} of the 2-sum. It is worth pointing out that this definition is slightly different from the usual definition of a 2-sum.

Notice that our definition of 2-sum can be easily extended to an operation between one graph $G_1$ and several other graphs $G_2^1$, $G_2^2$, ..., $G_2^t$, as long as each $G_2^i$ has a vertex of degree two and $G_1$ has $t$ pairwise nonadjacent vertices such that they all have degree two. The 2-{\it sum} of two classes ${\cal G}_1$ and ${\cal G}_2$ of graphs is the class $\mathcal G_1\oplus \mathcal G_2$ of all graphs obtained by 2-suming one graph from ${\cal G}_1$ with several graphs from ${\cal G}_2$ (so, in general, $\mathcal G_1\oplus \mathcal G_2\ne \mathcal G_2\oplus \mathcal G_1$). Now for any class ${\cal G}$ of graphs, let $\mathcal G^2 = \mathcal G \oplus \mathcal G$ and, inductively, $\mathcal G^{k+1} = \mathcal G^k \oplus \mathcal G$ for each integer $k\ge2$. Like before, a class of $k$-connected graphs is called {\it minor-closed} if every $k$-connected minor of a member is also a member. The following is the main result of this paper.

\begin{theorem}\label{thm:main}
Let ${\cal G}$ be a minor-closed class of $2$-connected graphs. Then ${\cal K}_2 \not \subseteq {\cal G}$ if and only if ${\cal G} \subseteq ( {\cal P} \cup {\cal A}_m )^m$ for some $m$.
\end{theorem} 

For graphs of higher connectivity, Theorem \ref{thm:main} can be simplified, as stated in the next corollary. For each positive integer $m$, let ${\cal A}_m'$ be the class of graphs obtained from graphs in ${\cal B}_m$ by adding disjoint strips so that their corners are identified with distinct vertices. 

\begin{corollary} 
Let ${\cal G}$ be a minor-closed class of $k$-connected graphs. \\ 
\indent{\rm(1)} If $k=3$, then ${\cal K}_2 \not \subseteq {\cal G}$ if and only if ${\cal G} \subseteq {\cal A}_m$ for some $m$;\\ 
\indent{\rm(2)} If $k=4$, then ${\cal K}_2 \not \subseteq {\cal G}$ if and only if ${\cal G} \subseteq {\cal A}_m'$ for some $m$; \\ 
\indent{\rm(3)} If $k=5$, then ${\cal K}_2 \not \subseteq {\cal G}$ if and only if ${\cal G} \subseteq {\cal B}_m$ for some $m$.
\end{corollary}

\begin{figure}[ht]
\centerline{\includegraphics[scale=1.0]{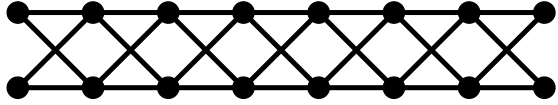}}
\caption{A regular strip}
\label{fig:strip}
\end{figure}

A strip is {\it regular} if it is obtained by 2-summing a sequence of $K_{3,3}\de e$ and then deleting the two vertices of degree two. Figure \ref{fig:strip} below illustrats such a graph. Since all fans and strips are minors of regular strips, we can use regular strips to describe all maximal 3-connected graphs that do not have a $K_{2,n}$-minor. For each pair of positive integers $a$ and $b$, let $J(a,b)$ be obtained by adding $a$ disjoint regular strips of length $b$ to $K_{4a}$.

\begin{corollary} 
There exists a function $a(n)$ such that every $3$-connected graph that does have a $K_{2,n}$-minor is a minor of $J(a(n),b)$ for some positive integer $b$.
\end{corollary}

The next corollary provides two sufficient conditions for having a large $K_{2,n}$-minor.

\begin{corollary} 
Let $n$ be a positive integer and let $G$ be a graph without $K_{2,n}$-minors. \\ 
\indent{\rm(1)} If $G$ is $3$-connected with minimum degree at least six then the size of $G$ is bounded.\\ 
\indent{\rm(2)} If $G$ is $4$-connected then the maximum degree of $G$ is bounded. 
\end{corollary}

We do not include in this paper the proofs of the last three corollaries since they are easy. The rest of the paper is devoted to proving Theorem \ref{thm:main}. First, in Section 2, we introduce internal 3-connectivity and we show how to decompose a 2-connected $K_{2,n}$-free graph into type-I graphs and internally 3-connected $K_{2,n}$-free graphs. Then, in Section 3, we present four lemmas, including a strengthening of Proposition \ref{prop:sn}. In sections 4, we show that, for any fixed integer $n\ge2$, every sufficiently large internally 3-connected $K_{2,n}$-free graph must contain a large fan or strip. In Section 5, we prove Theorem \ref{thm:main} for internally 3-connected graphs. Finally, we complete the whole proof in Section 6. 

We close this section by clarifying some terminology. For any integer $n\ge3$, let $C_n$ denote the cycle of length $n$. For any graph $G$, let $|G|$ denote its number of vertices and $||G||$ denote its number of edges. Let $d_G(x)$ denote the degree of a vertex $x$ of $G$. 
We will write $H\subseteq G$ if $H$ is a subgraph of $G$. Let $P\subseteq G$ be a path. If the two ends of $P$ are $x,y$ then we call $P$ an {\it $xy$-path}.  For any two vertices $u$ and $v$ of $P$, let $P[u,v]$ denote the unique $uv$-path contained in $P$. 
For any two graphs $G_1=(V_1,E_1)$ and $G_2=(V_2,E_2)$, let $G_1\cup G_2 =(V_1\cup V_2, E_1\cup E_2)$, $G_1\cap G_2 =(V_1\cap V_2, E_1\cap E_2)$, and $G_1\de G_2 = G_1\de V_2$. Finally, expression $(X\de Y)\de Z$ will be simplified as $X\de Y\de Z$.

\section{Decomposing 2-connected graphs}

A graph is called {\it internally $3$-connected} if it is obtained from a 3-connected graph by subdividing each edge at most once. 
Let $\cal C$ consist of $C_3$, $C_4$, $C_5$, $C_6$, $K_4\de e$, $K_{2,3}$, and internally 3-connected graphs. Let ${\cal C}_1=\{C_3,C_4,C_5,C_6,K_4\de e,K_{3,3}\de e\}$ and let ${\cal C}_2={\cal C}\de {\cal C}_1$. For each positive integer $n$, let ${\cal D}_n$ denote the class of all 2-connected graphs that do not have a $K_{2,n}$-minor, and let ${\cal D}_n' = {\cal D}_n\cap{\cal C}_2$. The main result in this section is the following lemma which asserts that graphs in ${\cal D}_n$ can be constructed from graphs in ${\cal D}_n'$ and $\cal P$ by taking 2-sums in at most $n$ iterations.

\begin{lemma} \label{lem:main2}
${\cal D}_n\subseteq ({\cal P}\cup {\cal D}_n')^n$, \ for all positive integers $n$.
\end{lemma}

We first show that every 2-connected graph can be expressed as 2-sums of graphs in $\cal C$ and then we prove the lemma using this decomposition. We proceed by proving a sequence of claims.

Two vertices $x$ and $y$ of a 2-connected graph $G$ form an {\it admissible $2$-cut} if $G$ can be expressed as the 2-sum of two graphs, for which both having fewer edges than $G$, and such that $\{x,y\}$ is the set of joins of the 2-sum.

\begin{lemma} \label{lem:21}
Let $G$ be $2$-connected and let $x,y\in V(G)$. Suppose $G\de\{x,y\}$ has $k\ge2$ components and $\{x,y\}$ is not admissible, then  \\ 
\indent {\rm (1)} $G=K_{2,3}$, if $k\ge3$;\\ 
\indent {\rm (2)} $G=K_4\de e$, if $k=2$ and $xy\in E(G)$;\\ 
\indent {\rm (3)} at least one component of $G\de\{x,y\}$ is an isolated vertex, if $k=2$ and $xy\not\in E(G)$.
\end{lemma}

\noindent{\it Proof}. Let $J_1$, $J_2$, ..., $J_k$ be all components of $G\de\{x,y\}$. Without loss of generality, let us assume that $|J_1| \le |J_2| \le ... \le |J_k|$. 

We first consider the case $k\ge3$. Let $H_1=G\de(J_3\cup ...\cup J_k)$, if $xy\not\in E(G)$, and $H_1=G\de e\de (J_3\cup ...\cup J_k)$, if $e=xy\in E(G)$. Let $H_2=G\de (J_1\cup J_2)$. For $i=1,2$, let $G_i$ be obtained from $H_i$ by adding a new vertex $z_i$ and two new edges $xz_i$ and $yz_i$. Clearly, $||G_2||<||G||$ and $G$ is the 2-sum of $G_1$ and $G_2$ with joins $x$ and $y$. Since $\{x,y\}$ is not admissible, we must have $||G_1||=||G||$. It follows that $k=3$, $|J_3|=1$, $xy\not\in E(G)$, and  thus $G=K_{2,3}$.

Next, we assume that $k=2$. Let $H_1=G\de J_2$; let $H_2=G\de J_1$, if $xy\not\in E(G)$, and $H_2=G\de e\de J_1$, if $e=xy\in E(G)$. For $i=1,2$, let $G_i$ be obtained from $H_i$ by adding a new vertex $z_i$ and two new edges $xz_i$ and $yz_i$. Clearly, $G$ is the 2-sum of $G_1$ and $G_2$ with joins $x$ and $y$. Since $\{x,y\}$ is not admissible, we must have $||G_i||=||G||$ for some $i$. It follows that either $|J_2|=1$, when $i=1$, or $|J_1|=1$ and $xy\not\in E(G)$, when $i=2$. Observe that $|J_1|=1$ holds in both cases, which implies (3). Finally, it is easy to see that $xy\in E(G)$ implies $i=1$ and thus $|J_2|=1$ and $G=K_4\de e$. \qed

A {\it labeled} graph consists of a 2-connected graph $G$ and a stable set $L$ for which all its vertices have degree two. An {\it admissible $2$-cut} of $(G,L)$ is an admissible $2$-cut $\{x,y\}$ of $G$ with $\{x,y\}\cap L=\emptyset$.

\begin{lemma} \label{lem:22} 
A $2$-connected graph $G$ belongs to $\cal C$ if and only if $(G,L)$ has no admissible $2$-cuts for some $L$.
\end{lemma}

\noindent{\it Proof}. Let $G\in \cal C$. If $G\ne C_6$, it is straightforward to verify that $(G,L)$ has no admissible 2-cuts for any $L$. If $G=C_6$, let $L$ be a stable set of size three. Then it is easy to see that $(G,L)$ has no admissible 2-cuts. Conversely, let $(G,L)$ have no admissible 2-cuts. By Lemma \ref{lem:21}, we may assume that, if $G\de\{x,y\}$ is disconnected for some $x,y\in V(G)\de L$, then $xy\not\in E(G)$, $G\de\{x,y\}$ has only two components, and exactly one of these two components is an isolated vertex. This assumption implies that, if $G$ has an $uv$-path $P$ such that $V(P)\de \{u,v\}\ne\emptyset$, both $u$ and $v$ have degree greater than two in $G$, and all other vertices of $P$ have degree two in $G$, then $G\de\{u,v\}$ has exactly two components, $P$ has exactly two edges, and $uv\not\in E(G)$. It follows that $G$ is either a cycle or an internally 3-connected graph. In the latter case, $G\in \cal C$. In the former case, let $G=C_k$ and let $x_1,x_2,...,x_k$ be the vertices of $G$ such that they are listed in the order they appear on the cycle. Without loss of generality, we may assume that $x_1\not\in L$. Let $t=\lceil k/2\rceil$. Then at least one of $x_t$ and $x_{t+1}$, say $x_t$, is not in $L$. If $k\ge7$, it is easy to check that $x_1$ and $x_t$ form an admissible 2-cut. Therefore, $k\le6$ and thus $G\in\cal C$. \qed

The next lemma is an obvious consequence of our definition of 2-sum. We omit its proof.

\begin{lemma}\label{lem:2sum}
If $G$ is a $2$-sum of $G_1$ and $G_2$ then both $G_1,G_2$ are minors of $G$.
\end{lemma}

We will need to talk about 2-summing more than two graphs. This sum can be obtained by taking 2-sums iteratively. In the following we make this more precise. 
As usual, if $\varphi$ is a function from a set $X$ to a set $Y$, then, for each $y\in Y$, we denote by $\varphi^{-1}(y)$ the set of elements $x\in X$ for which $\varphi(x)=y$. A {\it tree structure} $\Theta$ is a triple $(T,\ \{(G_t,L_t):t\in V(T)\},\ \varphi)$, where $T$ is a tree, each $(G_t,L_t)$ is a labeled graph, $\varphi$ is a function from a subset $dom(\varphi)$ of $\cup L_t$ to $E(T)$, and such that 

(i) $G_s$ and $G_t$ are disjoint if $s\ne t$;\\ 
\indent (ii) for each $e=st\in E(T)$, $\varphi^{-1}(e)$ has precisely two vertices with one in $L_s$ and one in $L_t$.

\noindent  For each $t\in V(T)$ and each $e\in E(T)$ that is incident with $t$, by (ii), $L_t$ has a unique vertex, which we denote by $l_{t,e}$, such that $\varphi(l_{t,e})=e$. Next, we explain why it makes sense to talk about the 2-sum of $\{G_t: t\in V(T)\}$ over the tree structure $\Theta$. 

Consider an edge $\alpha=uv\in E(T)$. Let us denote the new vertex of $T/\alpha$ by $\alpha$. Notice that, for each $x\in\{u,v\}$, the degree of $l_{x,\alpha}$ in $G_x$ is two, as $l_{x,\alpha}\in L_x$. It follows that the 2-sum of $G_u$ and $G_v$, over $l_{u,\alpha}$ and $l_{v,\alpha}$, is well defined. Let $G_{\alpha}$ be this 2-sum and let $L_{\alpha}=(L_u\cup L_v)\de\{l_{u,\alpha},l_{v,\alpha}\}$. Then it is easy to see that $L_{\alpha}$ is a stable set of $G_{\alpha}$ and all its vertices have degree two in $G_{\alpha}$. Therefore, $(G_{\alpha}, L_{\alpha})$ is a labeled graph. Let $\varphi/\alpha$ be the restriction of $\varphi$ to $dom(\varphi) \de \{ l_{u,\alpha}, l_{v,\alpha}\}$. Now it is straightforward to verify that $\Theta/\alpha = (T/\alpha, \ \{(G_t,L_t): t\in V(T/\alpha)\}, \ \varphi/\alpha )$ is a tree structure. 

Observe, from (i) above, that $(\Theta/\alpha)/\beta = (\Theta/\beta)/\alpha$, for any two distinct edges $\alpha$ and $\beta$ of $T$. Thus $\Theta/E(T)$ is well defined. Let $\Theta/E(T)=(K_1,\ \{(G,L)\},\ \varphi/E(T))$. Naturally, the labeled graph $(G,L)$ is called the {\it $2$-sum} of $\{(G_t,L_t): t\in V(T)\}$ over $\Theta$. Very often, we simply say that $(G,L)$ is the {\it $2$-sum} of $\Theta$. It follows from Lemma \ref{lem:2sum} that $G$ contains every $G_t$ as a minor. In fact, the following generalization is also true.

\begin{lemma} \label{lem:23}
Let $T'=(V',E')$ be a subtree of $T$, $\varphi'$ be the restriction of $\varphi$ to $\bigcup\{f^{-1}(e):e\in E'\}$, and $\Theta'=(T',\ \{(G_t,L_t):t\in V'\},\ \varphi')$. If the $2$-sum of $\Theta'$ is $(G',L')$, then $G'$ is a minor of $G$.
\end{lemma}

\noindent{\it Proof}. Let $T^*=T/E'$ and let its new vertex be $t'$. Let $(G_{t'},L_{t'})=(G',L')$; let $\varphi^*$ be the restriction of $\varphi$ to $dom(\varphi) \de dom(\varphi')$. Then $\Theta^*=(T^*,\ \{(G_t,L_t):t\in V(T^*)\},\ \varphi^*)$ is a tree structure and $G$ is the 2-sum of $\Theta^*$. By Lemma \ref{lem:2sum}, $G'=G_{t'}$ is a minor of $G$. \qed 

The next result, which is the first main part in proving Lemma \ref{lem:main2}, states that, roughly speaking, every 2-connected graph can be expressed as the 2-sum of graphs in $\cal C$ over a tree structure.

\begin{lemma} \label{lem:decomp} 
Every labeled graph $(G,L)$ is a $2$-sum of a tree structure $\Theta=(T$, $\{(G_t,L_t):$ $t\in V(T)\}$, $\varphi)$ such that each $G_t$ is in $\cal C$.
\end{lemma}

\noindent{\it Proof}. If there is a counterexample to the lemma then we can choose one with $||G||$ minimum. Clearly, $G$ is not in $\cal C$ because otherwise, $\Theta = (K_1,\{(G,L)\},\varphi)$ has the required property, where $\varphi$ is a function with an empty domain. Hence, by Lemma \ref{lem:22}, $(G,L)$ has an admissible 2-cut, say $\{x,y\}$. By definition, this means that $\{x,y\}$ is a subset of $V(G\de L)$ and $G$ is a 2-sum of two graphs $G^1$ and $G^2$, over $z^1\in V(G^1)$ and $z^2\in V(G^2)$, such that both $G^1$ and $G^2$ have fewer edges than $G$ and $\{x,y\}$ is the set of joins of the 2-sum. Let $L^1=(L\cup\{z^1\})\de V(G^2)$ and $L^2=(L\cup\{z^2\})\de V(G^1)$. Then it is easy to see that both $(G^1,L^1)$ and $(G^2,L^2)$ are labeled graphs. By the minimality of $G$, for each $i=1,2$, there exists a tree structure $\Theta^i=(T^i,\ \{(G_t,L_t): t\in V(T^i)\}, \ \varphi^i)$ such that $(G^i,L^i)$ is the 2-sum of $\Theta^i$ and each $G_t$ is in $\cal C$. Next, we construct from $\Theta^1$ and $\Theta^2$ a tree structure $\Theta$ and we prove that $\Theta$ has the required property. The existence of such a $\Theta$ implies that $G$ is not a counterexample and thus the lemma will be proved. 

For $i=1,2$, let $t^i$ be the vertex of $T^i$ for which $z^i\in V(G_{t^i})$. Let $T$ be obtained from $T^1$ and $T^2$ by adding an edge $t^1t^2$. Let $\varphi$ have domain $dom(\varphi^1) \cup dom(\varphi^2) \cup \{z^1,z^2\}$ such that $\varphi(z^1)=\varphi(z^2)=t^1t^2$ and $\varphi(z)=\varphi^i(z)$, where $i$ is the index for which $z\in dom(\varphi^i)$. Now it is routine to verify that $\Theta = (T,\ \{(G_t,L_t): t\in V(T)\}, \ \varphi)$ has the required property, as required. \qed

If the graph $G$ in Lemma \ref{lem:decomp} is in $\cal P$, then it is easy to show that each $G_t$ is in ${\cal C}_1$. We do not need this fact in our proof, but we do need the following.

\begin{lemma} \label{lem:25}
Let $\Theta = (T,\ \{(G_t,L_t):t\in V(T)\},\ \varphi)$ be a tree structure for which each $G_t$ is in ${\cal C}_1$. If $(G,L)$ is the $2$-sum of $\Theta$, then $G$ is in $\cal P$.
\end{lemma}

\noindent{\it Proof}. Clearly, each $G_t$ is in $\cal P$. In addition, it is straightforward to verify that the 2-sum of any two graphs in $\cal P$ is also in $\cal P$. Thus the result follows. \qed

Recall that ${\cal C}_2$ consists of $K_{2,3}$ and, except for $K_{3,3}\de e$, all internally 3-connected graphs. The next result says that graphs in ${\cal C}_2$ are useful in constructing a $K_{2,n}$-minor.

\begin{lemma} \label{lem:26}
Let $G \in {\cal C}_2$ and let $x$ and $y$ be two of its vertices of degree two. Then $G$ has two internally vertex-disjoint $xy$-paths $P$ and $Q$ such that $G\de\{x,y\}$ has a path $R$ on at least three vertices and is between $V(P)$ and $V(Q)$.
\end{lemma}

\noindent{\it Proof}. Let the neighbors of $x$ be $a$ and $b$ and let the neighbors of $y$ be $c$ and $d$. We first consider the case when $x$ and $y$ have a common neighbor, say $a=c$. Let $z$ be a neighbor of $a$ other than $x$ and $y$. Then it is easy to see that both $G\de a$ and $G\de\{a,z\}$ are connected. Thus $G\de\{a,z\}$ has an $xy$-path $P$ and $G\de a$ has a path $R'$ between $z$ and $V(P)$. Let $Q$ be the path on the three vertices $x$, $a$, $y$, and let 
$R$ be the path by adding the edge $az$ to $R'$. Then it is clear that the three paths $P$, $Q$, and $R$ have the required property. 

Next we assume that $x$ and $y$ do not have common neighbors. This implies that $G$ is internally 3-connected. We assume that the required paths do not exist and we deduce a contradiction by showing that $G=K_{3,3}\de e$. Clearly, there are two internally vertex-disjoint $xy$-paths, say, $P$ containing $a$ and $c$ and $Q$ containing $b$ and $d$. From the connectivity of $G$ we deduce that $x$ and $y$ are in the same component of $G\de\{b,c\}$ and thus there is an edge $pq$ with $p$ in $P\de \{x,y,c\}$ and $q$ in $Q\de \{x,y,b\}$. To simplify our argument, we choose such an edge with $P[x,p]$ as short as possible. Again, from the connectivity of $G$ we deduce that $x$ and $y$ are in the same component of $G\de\{p,q\}$ and thus either $G\de Q$ has a path between the two parts of $P\de\{x,y,p\}$, or $G\de P$ has a path between the two parts of $Q\de\{x,y,q\}$, or $G$ has an edge $p'q'$ with $p'$ on $P'$ and $q'$ on $Q'$, where $P'$ is the component of $P\de p$ that contains $c$ and $Q'$ is the component of $Q\de q$ that contains $b$. The first two cases do not occur since the three required paths do not exist. For the same reason, in the last case, we also have $pp'\in E(P)$ and $qq'\in E(Q)$. If $p\ne a$ or $q'\ne b$, then $G\de\{p,q'\}$ has an $xy$-path. It is routine to check that, no matter how this path goes, the three required paths always exist. Thus we must have $p=a$ and $q'=b$. Similarly, we also have $p'=c$ and $q=d$. Now the result follows from the connectivity of $G$. \qed

The following is an immediate consequence of Lemma \ref{lem:26} and Lemma \ref{lem:23}. 

\begin{lemma} \label{lem:27} 
Let $(G,L)$ be the $2$-sum of a tree structure $\Theta=(T,\ \{(G_t,L_t):t\in V(T)\},\ \varphi)$, where each $G_t$ is in $\cal C$. If some path of $T$ contains $n$ vertices $t$ with $G_t\in {\cal C}_2$, then $G$ has a $K_{2,n}$-minor.
\end{lemma}

\noindent{\it Proof of Lemma} \ref{lem:main2}. Let $G$ be a graph in ${\cal D}_n$. By Lemma \ref{lem:decomp}, $(G,\emptyset)$ is a $2$-sum of a tree structure $\Theta=(T,\ \{(G_t,L_t):t\in V(T)\},\ \varphi)$, where each $G_t$ is in $\cal C$. Let $U$ be the set of vertices $t\in V(T)$ such that $G_t\in {\cal C}_2$. Let $T_i=(V_i,E_i)$, $1\le i\le k$, be all the components of $T\de U$. Let $\Theta'=\Theta/(E_1 \cup E_2 \cup ... \cup E_k) = (T',\ \{(G_t,L_t):t\in V(T')\},\ \varphi')$ and let $U'=V(T')\de U$. Then we have the following observations. 

\begin{enumerate}[topsep=-2pt,itemsep=-2pt,leftmargin=40pt]

\item[(i)] By the definition of $T'$, distinct vertices of $U'$ are not adjacent in $T'$;

\item[(ii)]  By Lemma \ref{lem:25}, if $t\in U'$, then $G_t\in \cal P$; 

\item[(iii)]  By Lemma \ref{lem:23}, if $t\in U$, then $G_t\in {\cal D}_n'$;

\item[(iv)]  By Lemma \ref{lem:27}, each path of $T'$ has fewer than $n$ vertices in $U$. 
\end{enumerate}

\noindent It follows from (i) and (iv) that each path of $T'$ has fewer than $2n$ vertices. Therefore, by (ii) and (iii), $G$ belongs to $({\cal P}\cup{\cal D}_n')^n$, as required. \qed

\section{A few lemmas}

We present in this section four lemmas, which will be used in later sections. What these lemmas have in common is that they do not speak anything about $K_{2,n}$-minors.

Let $H$ be a subgraph of a graph $G$. An {\it $H$-bridge} of $G$ is a subgraph $B$ of $G$ such that either $B\cong K_2$ with $V(B)\subseteq V(H)$ yet $E(B)\not\subseteq E(H)$, or $B$ is formed by a component $C$ of $G\de H$ together with all edges between $C$ and $H$.  We will call $B$ {\it trivial} if it is the first kind. In case $E(H)=\emptyset$, an $H$-bridge could also be called a $V(H)$-bridge. The intersection of $V(B)$ and $V(H)$ is the set of {\it attachments} of $B$. 

Let $G$ be a graph and let $R\subseteq V(G)$. An $R$-{\it tree} is a tree $T\subseteq G$ such that all leaves of $T$ are in $R$. In case $T$ is a path we further require that none of its interior vertices is in $R$. Note that $G$ has an $R$-path (an $R$-tree that is a path) with ends $x,y$ if and only if $G$ has an $R$-bridge $B$ such that both $x,y$ are attachments of $B$. 
Our first lemma, a strengthening of Proposition \ref{prop:sn}, characterizes graphs that do not have an $R$-tree with many leaves.

\begin{lemma} \label{lem:rtree} 
There exists a function $f_{\ref{lem:rtree}}(n)$ with the following property. For any integer $n\ge2$, any connected graph $G$, and any $R\subseteq V(G)$, if $G$ does not have an $R$-tree with $n$ leaves, then $G$ has an induced subgraph $G_0$ such that 
\begin{enumerate}[topsep=-2pt,itemsep=-2pt]

\item[{\rm(1)}] $|R\cap V(G_0)|\le f_{\ref{lem:rtree}}(n)$;

\item[{\rm(2)}]  no two $G_0$-bridges have a common attachment; and  

\item[{\rm(3)}]  each $G_0$-bridge $B$ has precisely two attachments, both are in $R$, and every other vertex in $V(B)\cap R$ is a cut vertex of $B$ that separates the two attachments. 
\end{enumerate}
\end{lemma}

\noindent{\it Proof}. We prove that $f_{\ref{lem:rtree}}(n)=s^2$ satisfies the requirements, where $s=s(n)$ is the function determined in Proposition \ref{prop:sn}. To begin with, we define an auxiliary graph $H$ with vertex set $R$ and such that two vertices $x,y\in R$ are adjacent in $H$ if $G$ has an $R$-path $P_{x,y}$ with ends $x,y$. This definition implies immediately that $H$ is connected. It also follows that no tree $T\subseteq H$ may have $n$ leaves because otherwise the union of $P_{x,y}$, over all edges $xy$ of $T$, would contain an $R$-tree (of $G$) with $n$ leaves. By Proposition \ref{prop:sn}, $H$ is a subdivision of a graph $J$ with $|J|\le s$. 

If $H$ is not a cycle, let $R_0=\{x\in R: d_H(x)\ne 2\}$; if $H$ is a cycle, let $R_0$ consist of any single vertex of $J$. Let $R_1=\cup\{N_H[x]:x\in R_0\}$, where $N_H[x]$ consists of $x$ and all its neighbors in $H$. Since $R_0\subseteq V(J)$ and $||J||\le s(s-1)/2$, we have $|R_1|\le |J|+2||J||\le s^2$. Note that each $R_1$-bridge of $H$ is an $R_1$-path and no two nontrivial $R_1$-bridges of $H$ have a common attachment.

Let $\cal B$ consist of all $R$-bridges $B$ of $G$ for which all attachments of $B$ are contained in $R_1$. Note that if $B$ is an $R$-bridge of $G$ with at least three attachments then $B\in\mathcal B$ since these attachments form a clique in $H$.  
Let $G_0$ be the union of $B$ over all $B\in\mathcal B$. Then it is straightforward to verify that $G_0$ satisfies the requirements. \qed

The next is a Ramsey type result. For any two integers $x,y$, let $[x,y]$ denote the set of integer $z$ satisfying $x\le z\le y$. 

\begin{lemma}\label{lem:seq} 
Let $X$ and $Y$ be two finite sets of integers. If $|X|\ge n(m+1)$ then either \\ 
\indent $(1)$ there exist $x_1,...,x_n\in X$ and $y_1,...,y_n\in Y$ for which $x_1<y_1< ... <x_n<y_n$, or \\ 
\indent $(2)$ there exist  $x_1,...,x_m\in X$ for which $x_1<x_2<... <x_m$ and $Y\cap [x_1,x_m]=\emptyset$.
\end{lemma}

\noindent{\it Proof}. Let $k\ge0$ be the largest integer such that there exist $x_1,...,x_k\in X$ and $y_1,...,y_k\in Y$ for which $x_1<y_1< ... <x_k<y_k$. We assume without loss of generality that $k<n$ and the sequence $S=x_1,y_1,...,x_k,y_k$ is lexicographically minimum. Let $X_i=\{x\in X: x_i\le x\le y_i\}$ ($1\le i\le k$) and $X_{k+1}=\{x\in X: x>y_k\}$. Then the minimality of $S$ implies that $(X_1,...,X_{k+1})$ partitions $X$. 

Let $X_{k+1}=\{x_1',...,x_t'\}$, where $x_1'<x_2'<...<x_t'$. Then the minimality of $k$ implies that $[x_2',x_t']\cap Y=\emptyset$. Let us assume $t\le m$ because otherwise (2) holds. Thus $|X_i|\ge m+2$ holds for at least one $i\in\{1,...,k\}$. By the minimality of $S$ we have $[x_i,y_i]\cap Y \subseteq \{x_i,y_i\}$, and thus (2) holds since $(X_i\de\{x_i,y_i\})\cap Y=\emptyset$. \qed

A 0-1 matrix is {\it simple} if no two distinct columns are equal. The following, a special case of Lemma 3.2 of \cite{D}, determines all unavoidable large simple 0-1 matrices.

\begin{lemma} \label{lem:matrix} 
There exists a function $f_{\ref{lem:matrix}}(n)$ for which every simple 0-1 matrix with $f_{\ref{lem:matrix}}(n)$ columns can be permuted in such a way that the new matrix has an $n\times n$ submatrix $(a_{ij})$ such that \\ 
\indent {\rm(1)} $a_{ij}=1$ if and only if $i=j$; or \\ 
\indent {\rm(2)} $a_{ij}=1$ if and only if $i\ne j$; or \\ 
\indent {\rm(3)} $a_{ij}=1$ if and only if $i\ge j$.
\end{lemma}

The last lemma is known to many people, in various forms. 

\begin{lemma} \label{lem:local}
Let $G$ be a subdivision of a $3$-connected graph and let $x,y$ be distinct vertices of $G$. Then $G$ has an induced $xy$-path $P$ such that $G\de P$ is connected.
\end{lemma}

\noindent{\it Proof}. Let $P$ be an $xy$-path of $G$ such that the size of the largest $P$-bridge is maximized, where the size of a $P$-bridge $B$ is defined to be $|B\de P|$. Let $B_0$ be a bridge that attains this maximum. 
 
Let $\{x_1,....,x_k\}$ be the union of $\{x,y\}$ and the set of attachments of $B_0$. Suppose $x_1,...,x_k$ are listed in the order that they appear on $P$. 
For each $P$-bridge $B$, let $P(B)$ be the minimal path contained in $P$ that contains all attachments of $B$. Note that if $B$ is not $B_0$ then either $P(B)$ contains some $x_i$ as an interior vertex or $P(B)$ is a subgraph of some $P_i:=P[x_i,x_{i+1}]$. 

We claim that $G$ does not have a $P$-bridge $B\ne B_0$ such that some $x_i$ is an interior vertex of $P(B)$. Suppose $B$ is such a bridge. Let $u,v$ be the two ends of $P(B)$. Let $Q$ be a $uv$-path of $B$ that avoids all other attachments of $B$ and let $P'$ be obtained from $P$ by replacing $P[u,v]$ with $Q$. Note that $G$ has a $P'$-bridge $B'$ such that $V(B'\de P')\supseteq \{x_i\}\cup V(B_0\de P)$. This contradicts the choice of $P$ and thus the claim is prove. 

We assume that $G$ has a $P$-bridge $B_1\ne B_0$ because otherwise $P$ satisfies the conclusion of the lemma. The above claim implies $P(B_1)\subseteq P_i$ for some $i$. Let $G_1$ be the union of $P_i$ and all $P$-bridges $B$ with $P(B)\subseteq P_i$. Let $G_2$ be the graph formed by all edges of $G$ that are not in $G_1$. Then the above claim implies that $x_i,x_{i+1}$ are the only common vertices of $G_1$ and $G_2$. Since $G_1$ contains a cycle and since $G$ is a subdivision of a 3-connected graph, we deduce that $G_2$ is an $x_ix_{i+1}$-path. Therefore, $\{x_i,x_{i+1}\}=\{x_1,...,x_k\}=\{x,y\}$ and $B_0$ is a path satisfying the conclusion of the lemma. \qed

\section{Capturing large fans and strips} 

In this section we prove lemmas on producing a large fan or strip. The main result is Lemma \ref{lem:aug}, which says that every large internally 3-connected $K_{2,n}$-free graph must contain a big fan or strip.  

\begin{lemma} \label{lem:41}
Let $P$ be an induced path of $G$. Let $x\in V(G\de P)$ such that $G\de P\de x$ has $t$ components and $x$ has at least $(m+2)n^t$ neighbors in $P$. Then either $G$ has a $K_{2,n}$-minor or $P$ has a subpath $P^*$ such that $x$ has $m$ neighbors in $P^*$ and no other vertex outside $P$ has a neighbor in $P^*$.
\end{lemma}

\noindent {\it Proof}. Let $G_1,...,G_t$ be the components of $G\de P\de x$. Let $X$ be the set of vertices of $P$ that are adjacent to $x$; let $Y_i$ ($1\le i\le t$) be the set of vertices of $P$ that are adjacent to a vertex of $G_i$. Since $|X|\ge (m+1)n^t+n^{t-1}+...+n$, we deduce by repeatedly applying Lemma \ref{lem:seq} that either $G$ has a $K_{2,n}$-minor or $X$ contains two vertices $x_1,x_2$ such that $|Z\cap X|=m$ and $Z\cap Y_i=\emptyset$ ($1\le i\le t$), where $Z=V(P[x_1,x_2])$. The latter is exactly saying that $P^*:=P[x_1,x_2]$ satisfies the lemma. \qed 

For a path $P$, let $I(P)$ denote the set of interior vertices of $P$. Let $R(n_1,n_2,n_3)$ denote the least positive integer $n$ such that if edges of $K_n$ are colored by $1, 2, 3$, then there is an $i$-colored clique of size $n_i$ for some $i$.

\begin{lemma} \label{lem:strip} 
There exists a function $f_{\ref{lem:strip}}(n,m)$ with the following property. Let $P$ and $Q$ be disjoint induced paths of a graph $G$ for which $V(P\cup Q)=V(G)$. If $G$ has at least $f_{\ref{lem:strip}}(n,m)$ edges between $P$ and $Q$, then either $G$ has a $K_{2,n}$-minor or $P$ has a subpath $P^*$ and $Q$ has a subpath $Q^*$ such that \\ 
\indent {\rm(1)} \ vertices in $I(P^*)\cup I(Q^*)$ are only adjacent to vertices of $P^*\cup Q^*$; and \\ 
\indent {\rm(2)} \ $P^*$, $Q^*$, and edges between them form a fan or strip of length $\ge m$.
\end{lemma}

\noindent{\it Proof}. Let $f_{\ref{lem:strip}}(n,m)=R(r,r,(m+4)n^2)$, where $r=2^{f_{\ref{lem:matrix}}(2n)}s$ and $s=(m+6)n$. We prove that this function satisfies the lemma. 
We assume that $G$ does not have a $K_{2,n}$-minor and we prove the existence of $P^*$ and $Q^*$ that satisfy (1) and (2). We prove this by proving a sequence of claims. Let $F$ denote the set of edges between $P$ and $Q$. First, by Lemma \ref{lem:41}, we assume that

\noindent (i) {\it each vertex of $G$ is incident with fewer than  $(m+4)n^2$ edges of $F$.}

Let $V(P)=\{x_1,x_2,...,x_p\}$, $V(Q)=\{y_1,y_2,...,y_q\}$, and let the vertices be listed in the order that they appear on the paths. We will call two edges $x_iy_j$ and $x_{i'}y_{j'}$ {\it comparable} if $(i-i')(j-j')>0$. Let $K$ be the complete graph with vertex set $F$. For each $\alpha\in E(K)$, say $\alpha$ is between $x_iy_j$ and $x_{i'}y_{j'}$, we color $\alpha$ by $1$, $-1$, or $0$ if $(i-i')(j-j')$ is positive, negative, or zero, respectively. From (i) we know that $K$ does not have a clique of size $(m+4)n^2$ with all its edges colored 0. Therefore, $K$ has a clique of size $r$ for which its edges either are all colored by 1 or are all colored by $-1$. By reversing the ordering of vertices of $Q$, if necessary, we assume that $K$ has a clique of size $r$ with all its edges colored by 1. Equivalently,

\noindent (ii) {\it $F$ has a subset $M$ of $r$ pairwise comparable edges.}

We say that two edges $x_iy_j$ and $x_{i'}y_{j'}$ of $F$ {\it cross} if $(i-i')(j-j')<0$. We also say that the cross is {\it big} or {\it small} if $(i-i')(j-j')<-1$ or $(i-i')(j-j')=-1$, respectively. Now we define a graph $H$ with vertex set $F$ such that two members of $F$ are adjacent in $H$ if and only if they form a big cross. A subset $M'$ of $M$ is called {\it homogeneous} if, in $H$, every member of $F\de M'$ is adjacent to either all or none of members of $M'$. In the following, we prove that 

\noindent (iii) {\it $M$ contains a homogeneous set of size $s$.}

\begin{figure}[ht]
\centerline{\includegraphics[scale=0.7]{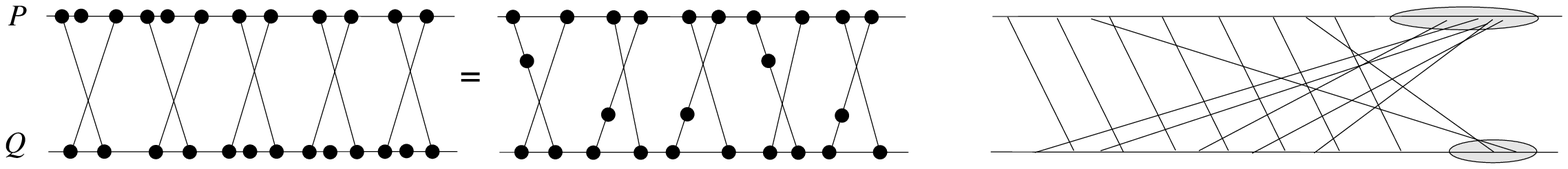}}
\caption{There exists a $K_{2,n}$-minor in both cases.}
\label{fig:chords}
\end{figure}

Let us call two vertices of $H\de M$ {\it similar} if they have the same set of neighbors in $M$. Clearly, being similar is an equivalence relation. Let $c$ be the number of equivalence classes. Then $M$ can be divided into $\le 2^c$ homogeneous sets. If $c< f_{\ref{lem:matrix}}(2n)$ then the largest homogeneous set satisfies (iii). Therefore, $c\ge f_{\ref{lem:matrix}}(2n)$ and thus $V(H)\de M$ contains a set $N$ of $f_{\ref{lem:matrix}}(2n)$ pairwise non-similar vertices. Let us define a 0-1 matrix $A=(a_{\alpha\beta})$ for which rows are indexed by members of $M$ and columns are indexed by members of $N$ and such that $a_{\alpha\beta}=1$ if and only if $\alpha$ and $\beta$ are adjacent in $H$. It follows from Lemma \ref{lem:matrix} that $A$ has a $2n\times 2n$ submatrix that satisfies one of (1-3) of Lemma \ref{lem:matrix}. Suppose rows and columns of this submatrix are indexed by $M'\subseteq M$ and $N'\subseteq N$. If (1) happens, then $P\cup Q$ together with members of $M'\cup N'$ form a subgraph of $G$ as illustrated on the left of Figure \ref{fig:chords}. In this case $G$ contains a $K_{2,2n}$-minor. Note that (2) does not happen since $M$ consists of pairwise comparable edges of $G$. If (3) happens, members of $N'$ can be naturally divided into two sets as illustrated on the right of Figure \ref{fig:chords}. In this case the bigger set together with $P\cup Q$ and $M'$ contains a $K_{2,n}$-minor. So (iii) is proved.

Let $M_1$ be a homogeneous set determined in (iii). Let $F_1=\{f\in F: f$ does not cross all $e\in M_1\}$. Since $M_1$ is homogeneous, if some $f\in F_1$ crosses some $e\in M_1$, then  $\{e,f\}$ must be a small cross. We will use this fact implicitly several time. Let $M_1'=\{e\in M_1: e$ does not cross any $f\in F_1\}$ and let $M_1''=\{e\in M_1: e$ crosses some $f\in F_1$ for which $f$ does not cross any $g\in F_1\de e\}$. Observe that each $e\in M_1''$ crosses a unique $f\in F_1$, and thus each edge in $F_1\de\{e,f\}$ crosses neither $e$ nor $f$. Let $M_1'''=M_1\de (M_1'\cup M_1'')$. Then each $e\in M_1'''$ crosses an edge $f_e\in F_1$ and $f_e$ crosses an edge $g_e\in F_1\de e$. Note that $\{f_e,g_e\}$ is a big cross. We show that 

\noindent (iv) {\it $|M_1'''|< 5n$.}

Suppose otherwise. Let $x_{i_1}y_{j_1}, ..., x_{i_{5n}}y_{j_{5n}}\in M_1'''$ with $i_1<...<i_{5n}$ and $j_1<...<j_{5n}$. Denote each $x_{i_k}y_{j_k}$ by $e_k$. Note that, for $k=3,4,...,5n-2$, the ends of $f_{e_k}$ are contained in $P[x_{i_{k-1}},x_{i_{k+1}}]\cup Q[y_{j_{k-1}},y_{j_{k+1}}]$ and the ends of $g_{e_k}$ are contained in $P[x_{i_{k-2}},x_{i_{k+2}}]\cup Q[y_{j_{k-2}},y_{j_{k+2}}]$. Therefore, $P\cup Q$ and $f_{e_{5k}},g_{e_{5k}}$ ($k=1,...,n$) form a subgraph of $G$ as illustrated on the left of Figure \ref{fig:chords}. This is impossible since $G$ would have a $K_{2,n}$-minor. Thus (iv) is proved.

Let $t=(m+1)n$. Since $s-5n=t$, we can choose $x_{i_1}y_{j_1}, ..., x_{i_{t+1}}y_{j_{t+1}}\in M_1'\cup M_1''$ with $i_1<...<i_{t+1}$ and $j_1<...<j_{t+1}$. Denote each $x_{i_k}y_{j_k}$ by $e_k$. For each $e_k\in M_1''$, let $f_k$ be the unique edge in $F_1$ that crosses $e_k$. Let $F_1'$ consist of all these $e_k$ and $f_k$. Then no edge in $F_1\de F_1'$ crosses any edge in $F_1'$. For $k=1,...,t$, let $Z_k$ be the set of edges $x_iy_j\in F_1\de F_1'$ such that $i_k\le i\le i_{k+1}$ and  $j_k\le j\le j_{k+1}$. If there are $n$ values of $k$ such that $Z_k$ contains a big cross, then these crossing edges and $P\cup Q$ form a subgraph of $G$ as illustrated on the left of Figure \ref{fig:chords}. This is impossible since $G$ would have a $K_{2,n}$-minor. Thus there are fewer than $n$ such $Z_k$. Since these $Z_k$ divide $e_1,...,e_{t+1}$ into at most $n$ consecutive sets, there must exist $k_0$ such that $Z_k$ does not contain a big cross for all $k=k_0,k_0+1, ..., k_0+m$. Let $P^*$ and $Q^*$ be the minimal paths in $P$ and $Q$, respectively, such that they contain all ends of $e_k,f_k$ for all $k=k_0, k_0+1, ..., k_0+m+1$. Then $P^*$ and $Q^*$ satisfy both (1) and (2) stated in Lemma \ref{lem:strip} and thus the lemma is proved. \qed

For any $H\subseteq G$, let $N(H)$ denote the set of vertices $v\in V(G\de H)$ such that $v$ is adjacent to at least one vertex of $H$. 

\begin{lemma}\label{lem:indp}
Let $G$ be internally 3-connected and $K_{2,n}$-free. Let $H\subseteq G$ such that both $H$ and $G\de H$ are connected. If $|N(H)| > (m+1)f_{\ref{lem:rtree}}(n)$ then $G\de H$ has a path $P$ of length at least $2m$ for which $d_{G\de H}(v)=2$ for all vertices $v$ of $P$.
\end{lemma}

\noindent{\it Proof}. Denote $R=N(H)$. Since $G$ is $K_{2,n}$-free, $G\de H$ does not contain an $R$-tree with $n$ leaves. By Lemma \ref{lem:rtree}, $G\de H$ has an induced subgraph $G_0$ that satisfies properties (1)-(3) of the lemma. It follows that $G\de H$ has at most $\frac{1}{2}f_{\ref{lem:rtree}}(n)$ $G_0$-bridges. Consequently, $G\de H$ has a $G_0$-bridges $B$ such that $B\de G_0$ contains at least $2m+1$ vertices of $R$. Since $G$ is internally 3-connected, we deduce from property (3) that $B$ is a path. Clearly, $P:=B\de G_0$ satisfies our requirements.  \qed

\begin{lemma}\label{lem:longp}
There exists a function $f_{\ref{lem:longp}}(n,m)$ such that every internally 3-connected $G$ with at least $f_{\ref{lem:longp}}(n,m)$ vertices must contain either a $K_{2,n}$-minor or an induced path $P$ of length $\ge m$ for which $G\de P$ is connected. 
\end{lemma}

\noindent {\it Proof}. Let $f_{\ref{lem:longp}}(n,m)=1+d+d^2+...+d^m$, where $d=(m+1)f_{\ref{lem:rtree}}(n)$. We prove that this function satisfies the lemma. 
If $G$ has a vertex $x$ of degree exceeding $d$ then the result follows from Lemma \ref{lem:indp} by letting $H$ be the single vertex $x$. So the maximum degree of $G$ is at most $d$. By our choice of $f_{\ref{lem:longp}}(n,m)$, $G$ has two vertices $x,y$ of distance $m$ between them. Then the result follows immediately from Lemma \ref{lem:local}. \qed

\begin{lemma}\label{lem:aug}
There exists a function $f_{\ref{lem:aug}}(n,m)$ such that every internally 3-connected $K_{2,n}$-free graph $G$ with at least $f_{\ref{lem:aug}}(n,m)$ vertices must contain two disjoint induced paths $P^*,Q^*$ that satisfy $(1)$ and $(2)$ of Lemma \ref{lem:strip}. 
\end{lemma}

\begin{figure}[ht] 
\centerline{\includegraphics[scale=0.6]{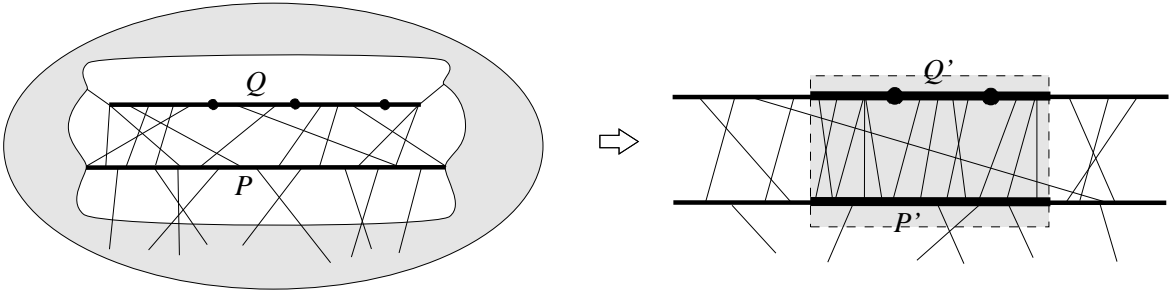}}
\caption{Paths $P$ and $Q$ contain paths $P'$ and $Q'$, respectively.}
\label{fig:findstrip}
\end{figure}

\noindent {\it Proof}. Let $f_{\ref{lem:aug}}(n,m)=f_{\ref{lem:longp}}(n,2(m+4)rn^r)$, where $r=(q+1) f_{\ref{lem:rtree}}(n)$, $q=f_{\ref{lem:strip}}(n,p)$, and $p=(m+6)n^4$.  
We prove that this function satisfies the lemma. By Lemma \ref{lem:longp}, $G$ has an induced path $P$ of length $2(m+4)rn^r$ such that $G\de P$ is connected. Let $R=N(P)$. If $|R|\le r$ then $R$ contains a vertex $x$ that has at least $(m+4)n^r$ neighbors in $P$. Note that each component of $G\de P\de x$ must contain a vertex from $R$, which implies that $G\de P\de x$ has fewer than $r$ components. Thus in the case $|R|\le r$ the result follows from Lemma \ref{lem:41}. Next, we assume $|R| > r$. By letting $H=P$ we deduce from Lemma \ref{lem:indp} that $G\de P$ has an induced path $Q$ of length $2q$ for which $d_{G\de P}(v)=2$ for all vertices $v$ of $Q$ (see Figure \ref{fig:findstrip}). Let $F$ be the set of edges between $P$ and $Q$. Then $|F|\ge q$. By Lemma \ref{lem:strip}, $P$ and $Q$ contain paths $P'$ and $Q'$, respectively, such that these two paths and edges between them form a fan or a strip of length at least $p$.

If $P'$ is a single vertex then $P^*:=P'$ and $Q^*:=Q'$ satisfy the lemma. Next we assume that $Q'$ is a single vertex $x$. Since $d_{G\de P}(x)=2$, $G\de P\de x$ has at most two components, which implies that $G\de P'\de x$ has at most four components. By Lemma \ref{lem:41}, $P'$ contains a path $P^*$ such that $x$ has $m+4$ neighbors in $P^*$ and no other vertex outside $P'$ has a neighbor in $P^*$. Thus $P^*$ and $Q^*:=Q'$ satisfy the lemma. Finally, we consider the case that neither $P'$ nor $Q'$ is a single vertex. 
Let $F'$ be a set of $p$ edges between $P',Q'$ such that they are pairwise noncrossing and nonadjacent.
Let $G'=G\de (F\de F')/E(Q')$ and let $x$ denote the resulting new vertex. By applying Lemma \ref{lem:41} to $G'$, $x$, and $P'$, we deduce that $P'$ contains a path $P''$ such that $x$ has $m+4$ neighbors in $P''$ and no other vertex in $G'\de P'$ has a neighbor in $P''$. Without loss of generality we assume that both ends $x_1,x_2$ of $P''$ are neighbors of $x$. Let $x_iy_i$ ($i=1,2$) be the unique edge in $F_1'$ that is incident with $x_i$ and let $Q''=Q'[y_1,y_2]$. Then it is straightforward to find $P^*\subseteq P''$ and $Q^*\subseteq Q''$ that satisfy the lemma. \qed

\section{Internally 3-connected graphs}

The purpose of this section is to prove the following.

\begin{theorem}\label{thm:i3c}
There exists a function $f_{\ref{thm:i3c}}(n)$ such that $\mathcal D_n'\subseteq \mathcal A_{f_{\ref{thm:i3c}}(n)}$ for all $n\ge4$.
\end{theorem}

We have established in Lemma \ref{lem:aug} that every large internally 3-connected $K_{2,n}$-free graph must contain a large fan or strip. To prove Theorem \ref{thm:i3c} it remains to show that maximal fans and strips are essentially disjoint. We prove this in a sequence of lemmas.

Let $G$ be obtained by adding a fan or a strip $J$ to a graph $H$. Let $Z$ be the set of corners of $J$. To simplify our notation, let us consider $J$ and $H$ as subgraphs of $G$. Note that $V(J\cap H)=Z$. Under these circumstances we will say that $J$ is a fan or strip of $G$. Vertices in $J\de Z$ are called {\it interior} vertices of $J$. Observe that for every interior vertex $v$ of $J$, all edges of $G$ that are incident with $v$ are contained in $J$. Also observe that $G$ may have edges with both ends in $Z$ yet not contained in $J$. A fan is called {\it nontrivial} if its length is at least two. 

\begin{lemma} \label{lem:fan1}
Let internally 3-connected graph $G$ contain two nontrivial fans $F_1,F_2$. If they have a common interior vertex $v$ then they must have the same center vertex.
\end{lemma} 

\noindent{\it Proof}. Suppose $z_1\ne z_2$ are centers of $F_1,F_2$, respectively. Let $P_i$ be the path $F_i\de z_i$. We first consider the case $d_G(v)=2$. Let $u,w$ be the two neighbors of $v$. Then each $P_i$ contains both $u,w$. Moreover, each $z_i$ is adjacent to both $u,w$. Since $F_1$ is nontrivial, at least one of $u,w$, say, $w$, is not a corner of $F_1$. Thus $w$ is an interior vertex of $F_1$, which implies $d_G(w)=3$ and so $wz_2\in E(P_1)$. Now from $uz_2\in E(G\de P_1)$ we deduce that both $u,z_2$ are corners of $F_1$, and thus they are the two ends of $P_1$. This is impossible since $F_1$ is nontrivial. Thus the case $d_G(v)=2$ is settled.

Now we assume $d_G(v)=3$. Then the set of neighbors of $v$ consists of $z_1,z_2$ and another vertex $z_3$. Note that each $P_i$ must contain $z_3$ and $z_{3-i}$. Since $d_G(z_{3-i})\ge4$, $z_{3-i}$ has to be a corner of $F_i$ and thus an end of $P_i$. It follows that $z_3$ is an interior vertex of both $F_1$ and $F_2$. From our discussion in the last paragraph we deduce that $d_G(z_3)\ne 2$. Thus $d_G(z_3)=3$ and $z_3$ is adjacent to $z_i$, the center of $F_i$, for both $i=1,2$. As a result, $z_3$ is an end of $P_i$, which means $P_i=z_3vz_{3-i}$, contradicting the nontriviality of $F_i$. \qed

Let $J$ be a fan of $G$. We call $J$ {\it maximal} if $G$ contains no other fan $J'$ with $J\subseteq J'$. 

\begin{lemma} \label{lem:2fan}
Let $F_1,F_2$ be distinct maximal nontrivial fans of an internally 3-connected graph $G$. If $v\in V(F_1\cap F_2)$ then $v$ is a corner of $F_i$ for both $i=1,2$, unless $G$ is obtained from a wheel by subdividing rims and at most one spoke.
\end{lemma}

\noindent{\it Proof}. Let us assume that $G$ is not one of the exceptions. Suppose $z\in V(G)$ such that $d_{G\de z}(v)=2$. We claim that $G\de z$ has a path $P$ containing $v$ and such that, for every $u\in V(P)$, $d_{G\de z}(u)=2$ if and only if $u$ is an interior vertex of $P$. To prove this claim we first observe that $G\de z$ has a path $P$ containing $v$ such that all its interior vertices $u$ satisfy $d_{G\de z}(u)=2$, because $v$ and its two neighboring edges (in $G\de z$) form such a path. Let us choose a maximal path $P$ with this property and let $x,y$ be the two ends of $P$. Clearly, we only need to show $d_{G\de z}(x)\ne 2\ne d_{G\de z}(y)$. Suppose otherwise. Without loss of generality, let $d_{G\de z}(x)=2$. Since $x$ is incident with exactly one edge in $P$, $G\de z$ must have an edge $xx'$ not in $P$. By the maximality of $P$, we must have $x'=y$. Since $G$ is not obtained from a wheel by subdividing rims, $G\de z$ must have a vertex outside $P$. It follows that $y$ is a cut vertex of $G\de z$. Since $G$ is internally 3-connected, $G\de P\de z$ can have only one vertex, which implies that $G$ is obtained from a wheel by subdividing rims and one spoke. This contradiction proves our claim.  

For $i=1,2$, let $z_i$ be the center of $F_i$, let $P_i$ be the path $F_i\de z_i$, and let $x_i,y_i$ be the two ends of $P_i$ (so $x_i,y_i,z_i$ are the corners of $F_i$). We first consider the case that $v$ is an interior vertex in both fans. In this case we deduce from Lemma \ref{lem:fan1} that $z_1=z_2$ and thus we can denote this vertex by $z$. Then we deduce from the above claim that $G\de z$ has a path $P$ such that $v$ is on $P$, $d_{G\de z}(u)=2$ for all internal vertices of $P$ and $d_{G\de z}(u)\ne2$ for the two ends of $P$. It follows that both $P_1,P_2$ are contained in $P$, which contradicts the maximality of $F_1,F_2$. 

It remains to consider the case that $v$ is corner of one fan but not a corner of the other. By symmetry we assume that $v$ is a corner of $F_1$ but is an interior of $F_2$. Since the degree (in $G$) of any interior vertex of a fan is at most three, $v$ cannot be $z_1$, and thus, by symmetry, we assume $v=x_1$. This implies $d_G(v)>2$ since $G$ is internally 3-connected. So $v$ has exactly three neighbors, one is $z_2$ and the other two are on $P_2$, which we call $u_1,u_2$. If both $vu_1,vu_2$ are in $F_1$ then one of $u_1,u_2$, say $u_1$, is $z_1$. Since $d_G(z_1)>3$, $u_1$ has to be a corner of $F_2$, which implies that $u_2$ is an interior vertex in both $F_1,F_2$. From our proof in the last paragraph we know that this is impossible, so some $vu_i$, say $vu_1$, is not contained in $F_1$. Since $F_1$ contains two edges incident with corner $v$, it follows that both $vu_2,vz_2$ are in $F_1$. Among the two neighbors of $v$ in $F_1$, one is a center and the other is an interior vertex, so one has degree exceeding three and the other has degree at most three. Therefore, $z_2=z_1$. Let $P$ be the path of $G\de z_i$ determined by the claim given in the beginning of this proof. Then both $P_1,P_2$ are contained in $P$, which contradicts the maximality of $F_1,F_2$. Thus $v$ must be a corner in both fans. \qed

\begin{figure}[ht] 
\centerline{\includegraphics[scale=0.7]{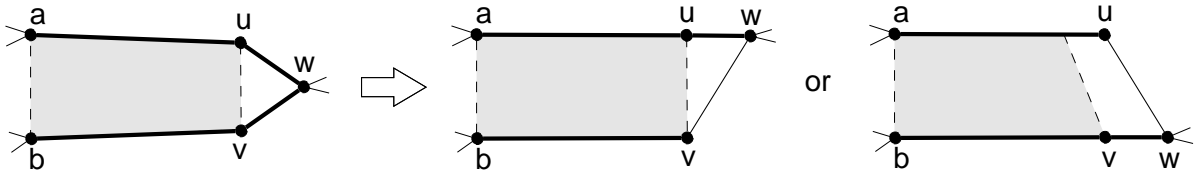}}
\caption{Two different ways to choose corners of a strip.}
\label{fig:boundary}
\end{figure}

Let $J$ be a strip of $G$. From the definition of trips we know that $J=H\de F$, where $H$ is a type-I graph with a reference cycle $C$ for which there are two nonadjacent edges $ab,cd$ of $C$ such that $F\subseteq \{ab,cd\}$. Let $P_1,P_2$ be the two paths of $C\de\{ab,cd\}$ and let $a,d$ be the two ends of $P_1$. In the rest of our discussions we need to specify the two paths and the four corners. We will write $(P_1,P_2;a,b,c,d)$ and call it the {\it boundary} of $J$. It should be pointed out that in some situations the choice of the four corners of $J$ is not unique. Consider strip $J$ illustrated in Figure \ref{fig:boundary}, where only three corners have neighbors outside $J$. In this scenario we may choose either $a,b,w,v$ or $a,b,w,u$ as the corners. In the following discussions, we insist that if $d_G(v)>d_G(u)=3$ then we require $v$ to be a corner, while if $d_G(v)>d_G(u)=2$ then we require $u$ to be a corner, and in the rest cases we may choose either $u$ or $v$ as the fourth corner. Note the the two choices do not affect the length of $J$.

Let $J'$ be a fan or strip of $G$. We say that $J'$ is {\it embedded } in $J$ if $J'\subseteq J$ and the following conditions are satisfied. If $J'$ is a fan with center $x$, we require that path $J'\de x$ is contained in one of $P_1,P_2$ and $x$ is contained in the other. If $J'$ is strip with boundary $(P_1',P_2';a',b',c',d')$ then we require that $P_i'\subseteq P_i$ ($i=1,2$), and $a'$ belongs to $P_1[a,d']$ if and only if $b'$ belongs to $P_2[b,c']$.

\begin{figure}[ht] 
\centerline{\includegraphics[scale=0.7]{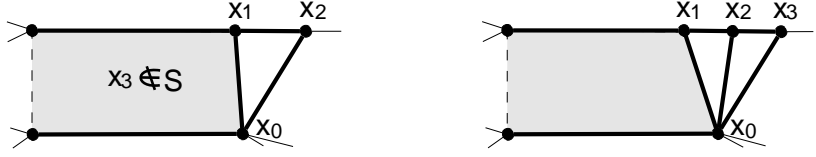}}
\caption{A basic fan $F$ agrees with a strip $S$.}
\label{fig:fanins}
\end{figure}

For $i=1,2$, let $J_i$ be a fan or a strip of $G$. We say that $J_1,J_2$ are {\it almost disjoint} if any vertex belongs to $J_1\cap J_2$ is a corner of both $J_1$ and $J_2$. 
We will call a fan {\it basic} if its length is two and it has exactly five vertices. If $F$ is a basic fan with center $x_0$ and $S$ is a strip, we say that $F$ {\it agrees} with $S$ if vertices of the path $F\de x_0$ can be expressed as $x_1x_2x_3x_4$ such that one of the two situations illustrated in Figure \ref{fig:fanins} occurs. In both cases, $x_0$ is a corner of $S$, and $x_1$ is in the interior of $S$ while $x_4$ is not. In the first case $x_2$ is a corner of $S$ and in the second case $x_3$ is a corner of $S$.  

\begin{lemma}\label{lem:smmetf}
Let $G$ contain a strip $S$ of length $\ge 6$ and a basic fan $F$. Then either $F$ is embedded in $S$, or $F,S$ are almost disjoint, or $F$ agrees with $S$. 
\end{lemma}

\noindent {\it Proof}. Let $F$ consist of a cycle $x_0x_1x_2x_3x_4x_0$ and two chords $x_0x_2,x_0x_3$. Let the boundary of $S$ be $(P_1,P_2;y_1,y_2,y_3,y_4)$. Then the length of $S$ implies immediately that

($*$) {\it any path of $S$ between $\{y_1,y_2\}$ and $\{y_3,y_4\}$ must have length $\ge3$}.

If $F$ is almost disjoint from $S$ then we are done. So we assume that some $x_i$ belongs to $S$ and this $x_i$ is not a corner of both $F$ and $S$. This assumption means $x_i$ is an interior vertex of at least one of $F$ and $S$. Consequently, it is straightforward to verify that at least one of $x_2,x_3$ is in $S$. 

We first consider the case that not both $x_2,x_3$ are in $S$. By symmetry we assume $x_3\not\in V(S)$. Then $x_2$ is a corner of $S$. Since $S$ contains at least two edges incident with $x_2$, $x_2x_0$ must belong to $S$. This implies that $x_0$ is a corner of $S$ since $x_0$ is adjacent to $x_3$. It also implies, by ($*$), that the set $\{x_0,x_2\}$ of corners must be either $\{y_1,y_2\}$ or $\{y_3,y_4\}$. Note that $x_4$ is not in the interior of $S$ since $x_4$ is adjacent to $x_3$. Thus $F$ agrees with $S$ as shown by the first graph in Figure \ref{fig:fanins}.

Now we assume both $x_2,x_3$ are in $S$. We prove that at least one of $x_2,x_3$ is in the interior of $S$. 
Suppose otherwise. Then both $x_2,x_3$ are corners of $S$. Note that either $x_2x_3\in E(S)$ or $x_2x_0$, $x_3x_0\in E(S)$. Thus, by ($*$), $\{x_2,x_3\} =\{y_1,y_2\}$ or $\{y_3,y_4\}$. We claim that $x_0$ is in the interior of $S$. If this is false, since $x_2$ is adjacent to an interior vertex of $S$, $x_1$ must be this vertex, so $x_0\in V(S)$, so $x_0$ is a corner of $S$, so path $x_0x_1x_2$ contradicts ($*$). By this claim we may assume $x_2=y_1$, $x_3=y_2$, and $x_0x_2\in E(P_1)$. It follows that $x_3x_4\in E(P_2)$ and then there is no room for $x_1$, which is a contradiction. So at least one of $x_2,x_3$ is in the interior of $S$.

By symmetry we assume that $x_2$ is an interior vertex of $P_1$. It follows that $x_0\in V(P_2)$ and thus path $x_1x_2x_3$ is contained in $P_1$. If $x_3$ is in the interior of $S$ then $F$ is embedded in $S$. So assume $x_3$ is a corner of $S$ and, by symmetry, let $x_3=y_4$. Since $S$ contains at least one edge of the form $y_3v$, where $v$ is on $P_1$, we deduce that $y_3=x_0$ or $x_4$. From the way we choose corners of a strip we must have $y_3=x_0$ (see Figure \ref{fig:boundary}). Moreover, since $x_4$ is adjacent to $x_3$, it is not in the interior of $S$. So $F$ agrees with $S$ as shown by the second graph in Figure \ref{fig:fanins}. \qed. 

\begin{lemma} \label{lem:domino}
Let $G$ have an induced subgraph $J$ consisting of a cycle $x_1x_2x_3x_4x_5x_6x_1$ and a chord $x_2x_5$. Suppose $J$ is a strip of $G$ with corners $x_1, x_3, x_4, x_6$. If $G$ contains a strip $S$ of length $\ge 8$ such that $x_2$ is an interior vertex of $S$, then $J$ is embedded in $S$. 
\end{lemma}

\noindent{\it Proof}. Let $(P_1,P_2;y_1,y_2,y_3,y_4)$ be the boundary of $S$. Then the length of $S$ implies

($*$) {\it any path of $S$ between $\{y_1,y_2\}$ and $\{y_3,y_4\}$ must have length $\ge4$}.

\noindent By symmetry we assume that $x_2$ is contained in $P_1$. We first prove that path $x_1x_2x_3$ is contained in $P_1$. Suppose this is false. We assume by symmetry that $P_1$ contains path $x_1x_2x_5$. We also assume that $x_1$ belongs to $P_1[y_1,x_2]$. If $x_6$ is not in $S$ then both $x_1,x_5$ are corners of $S$. It follows that $P_1$ has length two, contradicting ($*$). So $x_6$ belongs to $S$. In fact, by ($*$) again, $x_6$ must belong to $P_2$. By symmetry we assume that $x_6$ belongs to $P_2[y_2,x_3]$. 
Since $x_6x_5$ crosses $x_2x_3$ and since $x_6,x_3$ are not adjacent in $G$, $x_6x_5$ does not belong to $S$, which means that both $x_6,x_5$ are corners of $S$. Because of ($*$) and path $x_1x_2x_5$, $x_1$ must be an interior vertex of $P_1$ and thus $x_1x_6$ belongs to $S$. Therefore, path $x_6x_1x_2x_5$ contradicts ($*$). So path $x_1x_2x_3$ is contained in $P_1$. Consequently, $x_5$ belongs to $P_2$. We assume by symmetry that $x_1$ is contained in $P_1[y_1,x_2]$.

Our next goal is show that $x_5$ is an interior vertex of $P_2$. Suppose otherwise. By symmetry we assume $x_5=y_2$. Let $x_i$ be the unique neighbor of $x_5$ in $P_2$. There are two cases. Suppose $i=6$. Because of path $x_5x_2x_3$ we deduce from ($*$) that $x_3$ is an interior vertex of $S$, so $x_4$ must belong to $S$. Because of path $x_5x_2x_3x_4$ we deduce from ($*$) that $x_4\not\in\{y_3,y_4\}$. Again, since $x_3$ is an interior vertex of $S$, $x_4$ is not on $P_1[y_1,x_2]$. Consequently, $x_4$ is an interior vertex of $P_1[x_3,y_4]$ or $P_2[x_6,y_3]$. Since $x_4$ is adjacent to $x_5$, it does not belong to $P_2$. On the other hand, if $x_4$ belongs to $P_1$ then $x_1x_6$ would cross both $x_5x_2,x_5x_4$, a contradiction. So $i\ne6$. Now suppose $i=4$. Note that $S$ has no edge $y_1v$ with $v$ belongs to $P_2\de \{x_5,x_4\}$ because otherwise $y_1v$ would cross both $x_5x_2,x_4x_3$. In addition, $y_1$ is not adjacent to $x_4$ because otherwise $y_1$ would have to be $x_1$ (as $y_1x_4$ crosses $x_5x_2$), implying $x_1x_4\in E(G)$, which is not the case. So $y_1$ must be adjacent to $x_5$ since $y_1$ needs to have a neighbor in $P_2$. Therefore, $y_1=x_6$. This means that $x_1$ is an interior vertex of $S$ and thus all its neighbors are in $\{x_6,x_2\}\cup V(P_2)$. It follows that $d_G(x_1)=2$ because any extra edge would cross both  $x_5x_2,x_4x_3$. This contradicts the way we choose corners of a strip (see Figure \ref{fig:boundary}), which shows that $x_5$ is an interior vertex of $S$.

Now combining the last two paragraphs we deduce that paths $x_1x_2x_3$ and $x_4x_5x_6$ are contained in $P_1,P_2$ respectively. By symmetry we assume that $x_1$ belongs to $P_1[y_1,x_2]$. If $x_4$ belongs to $P_2[y_2,x_5]$, then at least one of $x_1x_6, x_3x_4$ is not in $S$ (since they both cross $x_2x_5$). It follows that the ends of this edge are corners of $S$, which contradicts ($*$). So $x_4$ belongs to $P_2[x_5,y_3]$. To finish showing that $J$ is embedded in $S$ we only need to prove, by symmetry, $x_1x_6\in E(S)$. This is clear if $x_1$ is an interior vertex of $S$. If $x_1$ is a corner of $S$, since $S$ should have an edge of the form $x_1v$ with $v\in V(P_2)$ and since this edge should not cross both $x_2x_5,x_3x_4$, $v$ must belong to $P_2[y_2,x_4]$. This implies $x_1x_6\in E(S)$ and thus our proof is complete. \qed

\begin{lemma} \label{lem:strip2}
Let $G$ contain a strip $J$ consisting of a cycle $x_1x_2...x_8x_1$ and two chords $x_2x_7,x_3x_6$, where $x_1,x_4,x_5,x_8$ are the corners. Suppose $J$ is not embedded in any strip of length $\ge6$. Then $J$ is almost disjoint from every strip $S$ of length $\ge8$. 
\end{lemma}

\noindent{\it Proof}. Suppose $S$ is a strip of length $\ge8$ such that some $x_i\in V(J\cap S)$ is not a corner of both $J$ and $S$. We only need to consider the case that $x_i$ is an interior vertex of $J$. This is because if $x_i$ is a corner of $J$ then it is an interior vertex of $S$. It follows that all neighbors of $x_i$, including an interior vertex $x_j$ of $J$, belong to $S$. So in this case we may take $x_j$, instead of $x_i$. 

Without loss of generality, we assume $x_2\in V(S)$. We prove that $J$ is embedded in a strip of length $\ge6$. Let $(P_1,P_2;y_1,y_2,y_3,y_4)$ be the boundary of $S$.

We first consider the case that $x_2$ is in the interior of $S$. By Lemma \ref{lem:domino}, $J\de\{4,5\}$ is embedded in $S$. We assume that paths $x_1x_2x_3$ and $x_8x_7x_6$ and contained in $P_1$ and $P_2$, respectively. Moreover, $x_1$ belongs to $P_1[y_1,x_2]$ and $x_8$ belongs to $P_2[y_2,x_7]$. If $x_3$ or $x_6$ is an interior vertex of $S$, then we deduce from Lemma \ref{lem:domino} that $J$ is embedded in $S$. So both $x_3,x_6$ are corners of $S$, which means $x_3=y_4$ and $x_6=y_3$. Note that none of $x_4,x_5$ is an interior vertex of $S$ because otherwise $x_3x_4$ or $x_5x_6$ could cross both $x_2x_7$ and $x_1x_8$, which is impossible. Now we choose interior vertices $y_1'$ and $y_2'$ from $P_1$ and $P_2$, respectively, such that $S$ contains a strip $S'$ with boundary $(P_1[y_1',y_4], P_2[y_2',y_3]; y_1',y_2',y_3,y_4)$ and such that the length of $S'$ is maximum. Then it is straightforward to verify that this length is $\ge5$. Let $S''$ be obtained from $S'$ by adding the path $x_3x_4x_5x_6$. Since $x_4,x_5\not\in V(S')$, it follows that $S''$ is a strip of length $\ge 6$. Moreover, $J$ is embedded in $S''$, as required.

Next we assume that no interior vertex of $J$ is in the interior of $S$. Thus $x_2$ is a corner of $S$. Let $x_2=y_1$. Since at least two of the three neighboring edges of $x_2$ are in $S$, at least one of $x_2x_3,x_2x_7$ is in $S$, which implies that one of $x_3,x_7$ is $y_2$. Since the neighbor of $x_2$ in $P_1$ is an interior vertex of $S$, this vertex must be $x_1$. So $x_8\in V(S)$. Since $x_8$ does not belong to $\{y_3,y_4\}$, as the path $x_2x_1x_8$ is too short, $x_8$ has to be an interior vertex of $S$ and thus $x_7$ belongs to $S$. But $x_7\not\in\{y_3,y_4\}$ because path $x_2x_1x_8x_7$ is still too short, we conclude that $x_7=y_2$. Consequently, $x_8$ is the unique neighbor of $x_7$ in $P_2$. If $x_3,x_6$ are not in $S$, then adding path $x_2x_3x_6x_7$ to $S$ results in a strip of length $\ge 9$ such that it contains $x_2$ as an interior vertex. Thus we deduce from the previous paragraph that the required strip exists. So we assume at least one of $x_3,x_6$ is in $S$. Then the first half of this paragraph implies that $\{x_3,x_6\}=\{y_3,y_4\}$ and $x_3x_4,x_5x_6$ are contained in $P_1\cup P_2$. Note that $V(J)\subseteq V(S)$ and $J\cup S$ is a component of $G$. Moreover, it is straightforward to verify that by deleting vertices and edges from the middle of $S$ we can find in $J\cup S$ a strip $S'$ of length $\ge6$ such that $J$ is embedded in $S'$. Now the proof is complete. \qed  

\begin{lemma}\label{lem:i3c}
Let $G$ be internally 3-connected and $K_{2,n}$-free and have no fans of length $\ge8$. Then $G$ is obtained from a graph on $<p$ vertices by adding $< np^2$ strips of length $\ge8$, where $p=f_{\ref{lem:aug}}(n,8)$.
\end{lemma}

\noindent{\it Proof}. We first prove that $G$ can be constructed from a graph on $<p$ vertices by adding strips of length $\ge8$. We prove this by induction on $|G|$. If $|G|<p$ then the result is trivially true. If $|G|\ge p$, by Lemma \ref{lem:aug}, $G$ contains a maximal strip $S_0$ of length $\ge8$. By deleting and contracting edges of $S_0$ we can reduce $S_0$ to a graph $J$ as describe in Lemma \ref{lem:strip2}, while maintaining the same corners. Let $G'$ be the resulting graph. Then $G'$ is internally 3-connected and $K_{2,n}$-free and have no fans of length $\ge8$. 
By our induction hypothesis, $G'$ can be constructed by adding strips $S_1, ..., S_t$ of length $\ge8$ to a graph $H$ with $|H|<p$. Since $S_0$ is maximal in $G$, $J$ is not embedded in any strip of $G'$ of length $\ge6$. So by Lemma \ref{lem:strip2}, each $S_i$ ($i\ge1$) is almost disjoint from $J$.
Thus $G'$ can be constructed by adding strips $J,S_1,...,S_t$ to $H'$, where $H'$ is obtained from $H$ by deleting the four interior vertices of $J$. Consequently, $G$ is obtained by adding strips $S_0,S_1,...,S_t$ to $H'$.

Now let $G$ be constructed from a graph $H$ with $<p$ vertices by adding strips $S_1,...,S_t$ of length $\ge8$. 
It remains to show that $t<np^2$. We define an auxiliary loopless graph $\Gamma$ as follows. Let $V(\Gamma)=V(H)$. For each $i$, we include six edges connecting all pairs of corners of $S_i$. So $\Gamma$ has exactly $6t$ edges. Since $G$ is $K_{2,n}$-free, no two vertices of $\Gamma$ are connected by $\ge n$ edges. Therefore, $t\le 6t\le (n-1) \binom{p-1}{2} < np^2$. \qed

\noindent{\it Proof of Theorem} \ref{thm:i3c}.  We prove that $f_{\ref{thm:i3c}}(n)=200 np^2$ satisfies the theorem, where $p=f_{\ref{lem:aug}}(n,8)$. 

Let $G\in\mathcal D_n'$. If $G= K_{2,3}$ then the result trivially holds since $f_{\ref{thm:i3c}}(n)>5$. So we assume $G$ is internally 3-connected. We also assume that $G$ is not obtained from a wheel by subdividing rims and at most one spoke, because otherwise either $|G|\le 10< f_{\ref{thm:i3c}}(n)$ or $G$ is constructed from a graph with at most seven vertices by adding a nontrivial fan. 

Let $\cal F$ be the set of maximal nontrivial fans of $G$. By Lemma \ref{lem:2fan}, any two members of $\cal F$ are almost disjoint. For every $F\in\cal F$, by contracting its rim edges we reduce $F$ to a basic fan $F'$. Let $G'$ be the resulting graph. It is clear that $G'$ is still internally 3-connected and $K_{2,n}$-free. Moreover, $G'$ contains no fans of length $\ge3$. By Lemma \ref{lem:i3c}, $G'$ is constructed from a graph $H$ on $<p$ vertices by adding trips $S_i$ ($i=1,...,t$) of length $\ell_i\ge8$, where $t< np^2$. 

It follows from Lemma \ref{lem:smmetf} that each $S_i$ contains a set $Z_i$ of at most four vertices that are from the two ends of $S_i$ such that $S_i\de Z_i$ is a strip and each $F'$ is either embedded in $S_i\de Z_i$ or almost disjoint from $S_i\de Z_i$. In fact, it is straightforward to verify that a similar $Z_i$ can be chosen so that \\ 
\indent (i) $Z_i\subseteq V(S_i)$ contains all corners of $S_i$, \\ 
\indent (ii) $S_i':=S_i\de Z_i$ is a strip embedded in $S_i$, \\ 
\indent (iii) each $F'$ is either embedded in $S_i'$ or almost disjoint from $S_i'$, \\ 
\indent (iv) the length of $S_i'$ is at least $\ell_i-6$ and $|Z_i|\le 100$. \\ 
We point out that a careful analysis can reduce the bound from 100 to 36, but we will not go through the trouble. The existence of $Z_1,...,Z_t$ implies that $G'$ can be constructed from a graph $H'$ on $\le |H|+96t<100np^2$ vertices by adding pairwise disjoint strips $S_1',...,S_t'$. Note that if we expand each basic fan embedded in $S_i'$ we obtain a strip $S_i''$ of $G$ with the same set of corners. 

Let $F_1,...,F_s\in \cal F$ such that $F_1',...,F_s'$ are all the basic fans of $G'$ that are almost disjoint from every $S_i'$. Since $F_1', ..., F_s'\subseteq H'$ and since interior vertices of $F_1', ..., F_s'$ are all distinct, so $s<|H'|/2$. For each $F_i$, let $u_i,v_i$ be its two noncenter corners. Then $F_i\de\{u_i,v_i\}$ contains a unique maximal fan $F_i''$. Note that $|F_i\de F_i''|\le 4$. So $G$ is obtained from a graph $H''$ on $\le |H'| + 2s < 200np^2$ vertices by adding fans $F_1'', ..., F_s''$ and strips $S_1'', ..., S_t''$. Note that non-center corners of the fans are all distinct and they are also different from centers of the fans and corners of the strips. Moreover, the corners of the strips are all distinct. Thus $G$ is an augmentation of $H''$. Now the proof is complete. \qed

\section{Proving Theorem \ref{thm:main}}

The main parts of the proof are done in previous sections. What we need to do in this section is to take care of the routine work that bridges the gaps. We first prove three lemmas concerning the ``if" part of Theorem \ref{thm:main}. 

\begin{lemma}\label{lem:typeI}
No graph in $\cal P$ contains a $K_{2,5}$-minor.
\end{lemma}

\noindent{\it Proof}. Suppose otherwise. Then we choose a counterexample $G$ with as few edges as possible. Since every outerplanar graph is $K_{2,3}$-free, $G$ must contain a pair $ab,cd$ of crossing chords. Observe that if $e$ is a crossing chord then both $G\de e$ and $G/e$ belong to $\cal P$. Thus the four edges on the cycle spanned by $ab$ and $cd$ are neither deleted nor contracted when the $K_{2,5}$-minor is produced. It follows that two of $a,b,c,d$, say $b$ and $d$, have degree two in the minor. Since $G\de d$ contains a $K_{2,4}$-minor, $G$ must contain another pair $a'b', c'd'$ of crossing chords avoiding $d$. Using the same argument we may assume that $d'$ is a degree two vertex in the minor. Now $G\de d\de d'$ contains a $K_{2,3}$-minor and thus $G$ must contain a third pair of crossing chords, which implies that the $K_{2,5}$-minor contains three edge-disjoint 4-cycles. This is a contradiction and thus the lemma is proved. \qed

\begin{lemma}\label{lem:sum}
If no graph in ${\cal G}_1$ has a $K_{2,m}$-minor, and no graph in ${\cal G}_2$ has a $K_{2,n}$-minor, then no graph in $\mathcal G_1\oplus \mathcal G_2$ has a $K_{2,mn}$-minor.
\end{lemma}

\noindent{\it Proof}. Let $G\in \mathcal G_1\oplus \mathcal G_2$ be obtained by 2-summing $G'\in\mathcal G_1$ with $G_1,...,G_t\in\mathcal G_2$, over $z_1',...,z_t'$ and $z_1,...,z_t$, respectively. Let the two neighbors of $z_i$ be $x_i,y_i$ and the two neighbors of $z_i'$ be $x_i',y_i'$; let $x_i^*,y_i^*$ be the identification of $x_i,y_i$ with $x_i',y_i'$, respectively. To simplify our notation, we assume $x_i=x_i'=x_i^*$ and $y_i=y_i'=y_i^*$. In other words, we consider $G_1\de z_1$, ..., $G_t\de z_t$, and $G'\de\{z_1',...,z_p'\}$ as subgraphs of $G$. Let $p$ be the largest integer such that $G$ has a $K_{2,p}$-minor. This minor can be expressed as two disjoint connected subgraphs $H_1,H_2$ of $G$ and distinct vertices $v_1, ..., v_p$ of $G\de (H_1\cup H_2)$ such that each $v_i$ has at least one neighbor in each of $H_1,H_2$. 

We first prove $V(H_1\cap G')\ne\emptyset \ne V(H_2\cap G')$. Suppose otherwise. Then we assume without loss of generality $H_1\subseteq G_1\de\{x_1,y_1,z_1\}$. It follows that $G_1\de z_1$ contains every $v_i$. If $H_2\subseteq G_1\de\{x_1,y_1,z_1\}$ then let $H_0:=H_2$; else let $H_0\subseteq G_1$ be obtained from $z_1$ and $H_2\cap G_1$ by adding $z_1u$ for every $u\in\{x_1,y_1\}\cap V(H_2)$. Note that $H_0$ is a connected subgraph of $G_1$ in both cases (to see this in the second case one may consider three subcases depending on $|\{x_1,y_1\}\cap V(H_2)|$). Thus  $H_0,H_1$ and $v_1,...,v_p$ create a $K_{2,p}$-minor in $G_1$, which implies $p<n\le mn$, as required.

Next we assume $V(H_1\cap G')\ne\emptyset \ne V(H_2\cap G')$. We classify $G_1,...,G_t$ into two groups, depending on if $G_i$ meets both $H_1$ and $H_2$. If some $G_i$, say $G_1$, meets at most one of $H_1$ and $H_2$, then $G_1$ does not contribute to the $K_{2,p}$-minor because the 2-sum of $G'$ with $G_2, ...,G_t$ already contains a $K_{2,p}$-minor. So we assume that each $G_i$ meets both $H_1$ and $H_2$, which implies that $x_i$ belongs to one of $H_1$ and $H_2$ and $y_i$ belongs to the other. As a result, both $H_1\cap G'$ and $H_2\cap G'$ are connected. Let $V_0=\{v_1,...,v_p\}\cap V(G')$ and $V_i=\{v_1,...,v_p\}\cap V(G_i)$ for $i=1,...,t$. Denote $p_i=|V_i|$. Then $H_1\cap G'$, $H_2\cap G'$, and $\{z_1',...,z_t'\}\cup V_0$ creates a $K_{2,t+p_0}$-minor of $G'$, which implies $t+p_0< m$. On the other hand, for each $i=1,...,t$, $H_1\cap G_i$, $H_2\cap G_i$, and $V_i$ creates a $K_{2,p_i}$-minor of $G_i$, which implies $p_i < n$. Now it is routine to verify $p=p_0+p_1+ ... + p_t\le (m-1)(n-1)<mn$. \qed

\begin{lemma}\label{lem:aug1}
There exists a function $f_{\ref{lem:aug1}}(n)$ such that every augmentation of a graph on at most $n$ vertices does not have a $K_{2,f_{\ref{lem:aug1}}(n)}$-minor.
\end{lemma}

\noindent{\it Proof}. By using the proofs from the last two lemmas one can obtain a linear function that satisfies the requirement. But we present a short proof for an exponential function $f_{\ref{lem:aug1}}(n) = 10n2^n$. Let $G$ be obtained from $G_0$ by adding strips $H_1,...,H_s$ and fans $H_{s+1}, ..., H_{s+t}$. Note that $4s + 2t\le |G_0|\le n$. Consider a strip $H_i$. Let $a,b,c,d$ be its four corners, where $ab$ and $cd$ belong to the reference cycle. Then $H_i$ has a set $A_i$ of at most four edges that separates $a,b$ from $c,d$. It follows that adding $H_i$ is equivalent to adding each part of $H_i\de A_i$ (which can be expressed as 2-sums) and then adding the edges of $A_i$. For each fan $H_i$, let $A_i$ consist of a rim edge of $H_i$. Then $A_i$ has the same property. In other words, $G\de A$, where $A=A_1\cup ... \cup A_{s+t}$, can be expressed as a 2-sum of $G_0'$ with $H_1',...,H_{2s+2t}'$, where $G_0'$ is obtained from $G_0$ by adding $2s+2t$ vertices of degree two, and $H_{2i-1}',H_{2i}'$ are obtained from the two parts of $H_i\de A_i$ by adding to each a degree two vertex. From the last two lemmas we deduce that $G\de A$ is $K_{2,5(n+2s+2t)}$-free. Since adding an edge may at most double the size of the larges $\mathcal K_2$-minor, $G$ does not have a $K_{2,m}$-minor, where $m\le 5(n+2s+2t)2^{|A|}\le f_{\ref{lem:aug1}}(n)$. \qed

Now we are ready to complete the proof of our main result.

\noindent{\it Proof of Theorem} \ref{thm:main}. The backward implication follows immediately from the last three lemmas, while the forward implication follows immediately from Lemma \ref{lem:main2} and Theorem \ref{thm:i3c}. \qed

\newpage

\noindent{\bf Acknowledgment.} This research is supported in part by NSF grant DMS-1500699.

\end{document}